\documentclass[titlepage,11pt]{article}
\usepackage{latexsym}
\usepackage{amsfonts}
\usepackage{amsmath}
\newcommand\blackslug{\hbox{\hskip 1pt \vrule width 4pt height 8pt depth 1.5pt
        \hskip 1pt}}
\newcommand\bbox{\hfill \quad \blackslug \bigbreak}
\def\DD{\hbox{-}}
\def\CC{\hbox{-}\cdots\hbox{-}}
\def\LL{,\ldots,}
\def\cupcup{\cup\cdots\cup}
\def\poly{poly-$\chi$-bounded}

\title{Polynomial bounds for chromatic number \\ VII. Disjoint holes}
\author{Maria Chudnovsky\thanks{Supported by  NSF grant DMS-2120644.}\\
Princeton University, Princeton, NJ 08544
\\
\\
Alex Scott\thanks{Research supported by EPSRC grant EP/V007327/1.}\\
Mathematical Institute, University of Oxford, Oxford OX2 6GG, UK
\\
\\
Paul Seymour\thanks{Supported by AFOSR grant
A9550-19-1-0187.}\\
Princeton University, Princeton, NJ 08544
\\
\\
Sophie Spirkl\thanks{We acknowledge the support of the Natural Sciences and Engineering Research
Council of Canada (NSERC), [funding reference number RGPIN-2020-03912].
Cette recherche a \'et\'e financ\'ee par le Conseil de recherches en sciences
naturelles et en g\'enie du Canada (CRSNG), [num\'ero de r\'ef\'erence
RGPIN-2020-03912].  }\\
University of Waterloo, Waterloo, Ontario N2L3G1, Canada}

\date{November 1, 2021; revised \today}

\newtheorem{thm}{}[section]

\newcommand{\Proof}{\noindent{\bf Proof.}\ \ }

\begin{document}
\maketitle
\begin{abstract}
A {\em hole} in a graph $G$ is an induced cycle of length at least four, and a {\em $k$-multihole} in $G$ is a set of pairwise disjoint and nonadjacent holes.  It is well known that if $G$ does not contain any holes then its chromatic number is equal to its clique number.  In this paper we show that, for any $k$, if $G$ does not contain a $k$-multihole, then its chromatic number is at most a polynomial function of its clique number.  We show that the same result holds if we ask for all the holes to be odd or of length four; and if we ask for the holes to be longer than any fixed constant or of length four.  This is part of a broader study of graph classes that are polynomially $\chi$-bounded.
\end{abstract}

\section{Introduction}
A function $\phi:\mathbb{N}\rightarrow\mathbb{N}$ is a {\em binding function} for a graph $G$ if $\chi(G)\le \phi(\omega(G))$,
where $\chi(G), \omega(G)$ denote the chromatic number of $G$ and the size of the largest clique of $G$, respectively.
A class $\mathcal{C}$ of graphs is {\em hereditary} if for every $G\in \mathcal{C}$, every graph isomorphic to an induced subgraph
of $G$ also belongs to $\mathcal{C}$. A hereditary class $\mathcal{C}$ is {\em $\chi$-bounded} if there is a function $\phi$
that is a binding function for each $G\in \mathcal{C}$,
and if so, we call $\phi$ a {\em binding function} for the class;
if there exists a polynomial binding function, we say that $\mathcal{C}$ is {\em \poly} (see \cite{survey} for a survey on $\chi$-bounded classes, and \cite{Schiermeyer} on \poly\ classes).  While many classes are known to be $\chi$-bounded, the proofs frequently give quite fast-growing functions, and it is natural to ask whether this is necessary.
A remarkable conjecture of Louis Esperet~\cite{esperet} asserted that every $\chi$-bounded hereditary class is \poly.
But this was recently disproved by 
Bria\'nski, Davies and Walczak \cite{BDW}.  So the question now is: which hereditary classes are \poly?

A hereditary graph class is defined by excluding some induced subgraphs. A graph is {\em $H$-free} if it has no induced subgraph isomorphic to $H$, and {\em $\{H_1,H_2\}$-free} means 
both $H_1$-free and $H_2$-free. 
There is a mass of results on $\chi$-bounded classes where one of the excluded graphs is a forest, but in this paper we 
consider some classes
where every excluded graph has a cycle. 
A {\em hole} is an induced cycle of length at least four, and {\em odd-hole-free} means containing 
no odd hole. A {\em four-hole} means a hole of length four.
%This remains open. For instance, the class of all $P_5$-free graphs is $\chi$-bounded, but we do not know that it is 
%\poly; and we even do not know that the class of all $\{P_5,C_5\}$-free graphs is \poly. 
Let us say a {\em $k$-multihole} of a graph $G$ is an induced subgraph 
with $k$ components, each a cycle of length at least four.
We denote the $k$-vertex path by $P_k$ and the $k$-vertex cycle by $C_k$. 

Graphs with no $1$-multihole are chordal and hence perfect.
The class of graphs with no $k$-multihole in which all the cycles have odd length, is shown in~\cite{ind7}
to be $\chi$-bounded, but it contains the class of $\{P_5,C_5\}$-free graphs, and we cannot yet prove it is \poly\  
(see \cite{poly4} for the best current bounds).
If we replace ``odd'' by ``long'', the same applies: it is shown in~\cite{ind10} that for every $\ell\ge 0$, 
the class of graphs with no $k$-multihole in which all the cycles have length at least $\ell$ is
$\chi$-bounded (and we cannot yet prove it is \poly, for the same reason).
But we can if we permit cycles of length four to be components of the multiholes we are excluding. 
We will show:
\begin{thm}\label{oddandfourmultihole}
	For each integer $k\ge 0$, let $\mathcal{C}$ be the class of all graphs $G$ with 
	no $k$-multihole in which every component either has length four or odd length. Then $\mathcal{C}$ is \poly.
\end{thm}

If we change ``odd'' to ``long'', it also works:
\begin{thm}\label{longandfourmultihole}
        For all integer $k\ge 0$ and $\ell\ge 4$, let $\mathcal{C}$ be the class of all graphs $G$ with 
        no $k$-multihole in which every component either has length four or length at least $\ell$. Then $\mathcal{C}$ is \poly.
\end{thm}

This second one we can make stronger (we could not prove
the corresponding strengthening
of the first):
\begin{thm}\label{longandbip}
        For all integers $k, s\ge 0$, and $\ell\ge 4$, let $\mathcal{C}$ be the class of all graphs $G$ such that no induced subgraph of $G$
	has exactly $k$ components, each of which is either isomorphic to $K_{s,s}$ or a cycle of length at least $\ell$.
	Then $\mathcal{C}$ is \poly.
\end{thm}
(In general, $K_{s,t}$ denotes the complete bipartite graph with parts of cardinality $s$ and $t$.) 
Both these results derive from a theorem about $K_{s,s}$, which we will explain 
in the next section.

\section{Excluding a disjoint union, and self-isolation}

If $A\subseteq V(G)$, $G[A]$ denotes the subgraph of $G$ induced on $A$; and
we write $\chi(A)$ for $\chi(G[A])$ and $\omega(A)$ for $\omega(G[A])$. Two disjoint subsets of $V(G)$ are {\em anticomplete} 
if there are no edges
between them, and {\em complete} if every vertex of the first subset is adjacent to every vertex of the second. A graph $G$ {\em contains} a graph $H$ if some
induced subgraph
of $G$ is isomorphic to $H$, and such a subgraph is a {\em copy} of $H$.
A function $\phi:\mathbb{N}\rightarrow\mathbb{N}$ is {\em non-decreasing} if $\phi(x)\le \phi(y)$ for all $x,y\in \mathbb{N}$ with $x\le y$.

Let us say a graph $H$ is {\em self-isolating} if for every non-decreasing
polynomial $\psi:\mathbb{N}\rightarrow \mathbb{N}$, there is a polynomial
$\phi:\mathbb{N}\rightarrow \mathbb{N}$ with the following property. For every graph $G$ with $\chi(G)> \phi(\omega(G))$,
there exists $A\subseteq V(G)$ with $\chi(A)>\psi(\omega(A))$, such that either
\begin{itemize}
	\item $G[A]$ is $H$-free, or
	\item $G$ contains a copy $H'$ of $H$ such that $V(H')$ is disjoint from and anticomplete to $A$.
\end{itemize}

Self-isolation is of interest in considering polynomial $\chi$-boundedness for the class of $H$-free graphs, where $H$ is a forest. Say a forest 
$H$ is {\em good} if the class of $H$-free graphs is polynomially $\chi$-bounded.  It might be
true that every forest is good (strengthening the Gy\'arf\'as-Sumner conjecture~\cite{gyarfas,sumner} from $\chi$-boundedness to polynomial $\chi$-boundedness),
but this has only been proved for a few simple kinds of tree $H$,
and some (not all) of the forests that are disjoint unions of these trees. It is not known that if trees $H_1,H_2$ are good, then the disjoint
union of $H_1$ and $H_2$ is good. For instance,
trees of diameter three are good~\cite{poly3}, but disjoint unions of them might not be as far as we know.
But self-isolation helps here: if $H_1$ and $H_2$ are good forests, and one of them is self-isolating, then the disjoint union 
of $H_1$ and $H_2$ is good. Some good trees are known to be self-isolating (namely, stars and four-vertex paths), so we can happily 
take disjoint unions with them and preserve goodness.

Which graphs are self-isolating? We know very little at the moment: there are very few graphs that we know to have the property, and none
that we know not to have the property. 
(Could it be that all graphs are self-isolating? Certainly, if we change the definition of
self-isolating, replacing the polynomials $\phi,\psi$ by general functions, it is easy to show that all graphs have the property, by induction
on $\omega(G)$.)
A graph is self-isolating if all its components 
are self-isolating, but the only connected graphs that we know are self-isolating are complete graphs (proved below), paths 
of arbitrary length (proved in~\cite{poly6}), and  
complete bipartite graphs (proved in the next section).
The main result of~\cite{poly2} was that stars are self-isolating, so our result that complete bipartite graphs
are self-isolating generalizes this. 
The last takes up the main part of this paper, and is most of
what we need to prove
\ref{oddandfourmultihole} and \ref{longandbip}. 

First, complete isolation:
\begin{thm}\label{completeisolation}
	Every complete graph is self-isolating.
\end{thm}
\Proof (This proof was derived from a similar proof in~\cite{Schiermeyer2}.) Let  $\psi:\mathbb{N}\rightarrow \mathbb{N}$  be a 
non-decreasing polynomial, and let $H$ be a $k$-vertex complete graph. Let $\phi$ be the polynomial $\phi(x)=(x+1)^{k}\psi(x)+x$ for $x\in \mathbb{N}$.
Now let $G$ be a graph with chromatic number more than 
$\phi(\omega(G))$, and let $K$ be a clique of $G$ with cardinality $\omega(G)$. If $\omega(G)<k$, then the first bullet in the definition
of self-isolating holds, so we assume that $\omega(G)\ge k$. For each $X\subseteq K$ with $|X|=k$, let $A_X$ be the set of vertices in
$V(G)\setminus K$ that are nonadjacent to every vertex in $X$; and for every $Y\subseteq K$ with $|Y|=k-1$, let $B_Y$
be the set of vertices in $V(G)\setminus K$ that are adjacent to every vertex in $K\setminus Y$. Thus $V(G)\setminus K$ is the union
of the $\binom{\omega(G)}{k}$ sets $A_X$ and the $\binom{\omega(G)}{k-1}$ sets $B_Y$; and since 
$$\binom{\omega(G)}{k}+\binom{\omega(G)}{k-1}=\binom{\omega(G)+1}{k}\le (\omega(G)+1)^{k},$$
and $\chi(G\setminus K)>(\omega(G)+1)^{k}\psi(\omega(G))$, 
one of the sets $A_X$ or $B_Y$ has chromatic number more than $\psi(\omega(G))$. If $\chi(A_X)>\psi(\omega(G))$ for some $X$, 
then $G[X]$
is a copy of $H$ anticomplete to $A_X$, and since $\psi(\omega(G))\ge \psi(\omega(A_X))$, the second bullet in the definition 
of self-isolating holds. If $\chi(B_Y)>\psi(\omega(G))$
for some $Y$, then since $|K\setminus Y|=\omega(G)-k+1$ and $B_Y$ is complete to $K\setminus Y$, it follows that
$\omega(B_Y)<k$ and so $G[B_Y]$ is $H$-free, and the first bullet in the definition of self-isolating holds.
This proves \ref{completeisolation}.~\bbox

\section{Complete bipartite isolation}

We turn to the proof that 
\begin{thm}\label{bipisol}
	Every complete bipartite graph is self-isolating.
\end{thm}
We will in fact prove something a little stronger.
Let $\psi:\mathbb{N}\rightarrow\mathbb{N}$ be some non-decreasing function.
An induced subgraph $H$ of a graph $G$ is {\em $\psi$-nondominating} if there exists a set $A\subseteq V(G)$ disjoint from and 
anticomplete to $V(H)$, with $\chi(A))\ge \psi(\omega(A))$.
If $\psi:\mathbb{N}\rightarrow\mathbb{N}$ is a non-decreasing function and $q\ge 0$ is an integer,
a {\em $(\psi, q)$-sprinkling} in a graph $G$ is a pair $(P,Q)$ of disjoint subsets of $V(G)$, such that
\begin{itemize}
	\item $\chi(P)> \psi(\omega(P))$; and
	\item $\chi(Q)> \psi(\omega(Q))+qr$, where $r$ is the maximum over $v\in P$ of the chromatic number of the set of 
		neighbours of $v$ in $Q$.
\end{itemize}
(This is closely related to what was called a ``$(\psi, q)$-scattering'' in~\cite{poly6}.)
We will prove:
\begin{thm}\label{strongisol}
	Let $s,q\ge 0$ be integers, and let $\psi:\mathbb{N}\rightarrow\mathbb{N}$ be a non-decreasing polynomial.
	Then there is a polynomial $\phi:\mathbb{N}\rightarrow\mathbb{N}$ with the following property. For every graphs $G$ with 
with $\chi(G)>\phi(\omega(G))$, either:
	\begin{itemize}
		\item there is a $\psi$-nondominating copy of $K_{s,s}$ in $G$, or
		\item there is a  $(\psi, q)$-sprinkling in $G$.
	\end{itemize}
\end{thm}
\noindent{\bf Proof of \ref{bipisol}, assuming \ref{strongisol}.\ \ }
Let $s,s'\ge 0$ be integers, where $s'\le s$. We will show that $K_{s,s'}$ is self-isolating. (It is not enough to show this
when $s=s'$, because we  do not know that every induced subgraph of a self-isolating graph is self-isolating.)
Let $\psi:\mathbb{N}\rightarrow\mathbb{N}$ be a non-decreasing polynomial, let
$q=s+s'$, and let $\phi$ satisfy \ref{strongisol}. Let $G$
        be a graph with $\chi(G)>\phi(\omega(G))$. We claim that either there is a $\psi$-nondominating copy of $K_{s,s'}$ in $G$, or
there exists $A\subseteq V(G)$ with $\chi(A)>\psi(\omega(A))$ such that $G[A]$ is $K_{s,s'}$-free.
If there is a $\psi$-nondominating copy of $K_{s,s}$ in $G$, then there is also one of $K_{s,s'}$, so by \ref{strongisol}, 
we may assume that there is a  $(\psi, q)$-sprinkling $(P,Q)$ in $G$. If $G[P]$ is $K_{s,s'}$-free, the claim holds,
so we assume that there is a copy $H$ of $K_{s,s'}$ in $G[P]$. Thus $|H|=q$. 
Let $r$ be the maximum over $v\in P$ of the chromatic number of the 
set of neighbours of $v$ in $Q$. 
The set of vertices in $Q$ with a neighbour in $V(H)$ has chromatic number at most $|H|r=qr$; and $\chi(Q)> \psi(\omega(Q))+qr$
from the definition of a $(\psi, q)$-sprinkling. Consequently $H$ is $\psi$-nondominating, and hence $K_{s,s'}$ is self-isolating.~\bbox

To prove \ref{strongisol} we will need the following lemma:
\begin{thm}\label{bignonnbr}
        For every graph $G$ that is not a complete graph, there is a vertex $v$ such that the set of vertices different
        from and nonadjacent to
        $v$ has chromatic number at least $\chi(G)/\omega(G)$.
\end{thm}
\Proof
Let $X$ be a maximum clique of $G$, and for each $x\in X$, let $D_x$ be the set of vertices of $G$ different from
and nonadjacent to $x$. Since $G$ is nonnull, it follows that $X\ne \emptyset$. But $V(G)$ is the union of the sets
$D_x\cup \{x\}$ over $x\in X$, because of the maximality of $X$;
and so there exists $v\in X$ such that $\chi(D_v\cup \{v\})\ge \chi(G)/\omega(G)$. Choose such
a vertex $v$ with $D_v\ne \emptyset$ if possible. If $D_v\ne \emptyset$, then $\chi(D_v\cup \{v\})=\chi(D_v)$,
since there are no edges between $v$ and $D_v$, and so the theorem holds. Thus we may assume (for a contradiction)
that $D_v=\emptyset$, and
so $1=\chi(D_v\cup \{v\})\ge \chi(G)/\omega(G)$. Since $\chi(G)/\omega(G)\ge 1$, equality holds, and
so $\chi(D_x\cup \{x\})\ge \chi(G)/\omega(G)$ for every $x\in X$; and so $D_x=\emptyset$ for all $x\in X$, from the
choice of $v$. Consequently $V(G)=X$, and $G$ is a complete graph, a contradiction. This proves \ref{bignonnbr}.~\bbox

The proof of \ref{strongisol} will be by examining the largest ``template'' in $G$. With $s$ fixed, let us say that, for all 
integers $t,k\ge 0$,  a {\em $(t,k)$-template} 
in $G$ is a sequence
$(A_1\LL A_k)$ of pairwise disjoint subsets of $V(G)$, each of cardinality $t$, such that for $1\le i<j\le k$, and 
for every stable set $S\subseteq A_j$ with $|S|=s$,
every vertex in $A_i$ has a neighbour in $S$. The next result will enable us to find a $(t,2)$-template.
If $v\in V(G)$, we denote the set of neighbours of a vertex $v$ by $N(v)$ or $N_G(v)$.

\begin{thm}\label{bigbip}
	Let $s,q,t\ge 0$ be integers, and let $\psi:\mathbb{N}\rightarrow\mathbb{N}$ be a non-decreasing polynomial.
	Let $G$ be a graph with 
	\begin{align*}
		\chi(G)&> \omega(G)^{s}\left(\left(s+t^s\right)\psi(\omega(G))+t\right) \text{ and }\\
	%$$\chi(G)> \omega(G)^{s}(t^sp(\omega(G))+t)+ (1+\omega(G)+\cdots+\omega(G)^{s-1})\omega(G)p(\omega(G))$$
		\chi(G)&\ge q^{s} t+ \left(2+q+q^2+\cdots+q^{s-1}\right)\psi(\omega(G)) +2.
	\end{align*}
	Then either
        \begin{itemize}
                \item there is a $\psi$-nondominating copy of $K_{s,s}$ in $G$, or
                \item there is a  $(\psi, q)$-sprinkling in $G$, or
		\item $G$ contains a $(t,2)$-template.
        \end{itemize}
\end{thm}
\Proof
We may assume that $s,t\ge 1$.
Define $p=\psi(\omega(G))$. For $0\le i\le s$, define 
\begin{align*}
	m_{i}&= \omega(G)^{s-i}\left(t^sp+t\right)+ \left(1+\omega(G)+\cdots+\omega(G)^{s-i-1}\right)p\\
	n_{i}&=q^{s-i} t+ \left(1+q+q^2+\cdots+q^{s-i-1}\right)p.
\end{align*}
Thus $m_s=t^sp+t$, and $m_{i}=\omega(G)m_{i+1}+ p$ for $0\le i<s$; and
$n_s=t$ and $n_i=qn_{i+1}+p$ for $0\le i<s$.
By hypothesis, $\chi(G)>m_0$ and $\chi(G)> n_0+p+1$.
\\
\\
(1) {\em There is a vertex $v_1$ such that $\chi(N(v_1))> n_1$ and $\chi(M(v_1))>m_1$, where 
$M(v_1)=V(G)\setminus (N(v_1)\cup \{v_1\})$.}
\\
\\
Let $S$ be the set of all vertices $v$ with
$\chi(N(v))\le n_1$. If $\chi(S)>p$, choose a subset $P\subseteq S$ with $\chi(P)= p+1$, and let $Q=V(G)\setminus P$. Then 
$$\chi(Q)\ge \chi(G)-(p+1)> n_0=p+qn_1,$$ 
and so $(P,Q)$ is a $(\psi,q)$-sprinkling. We therefore assume that $\chi(S)\le p$.
%\chi(G)\ge 2p_0+qn_1+2$,
Let $R=V(G)\setminus S$. Thus 
$$\chi(R)\ge \chi(G)-p>m_0-p=\omega(G)m_{1}\ge \omega(G),$$ 
and so $R$ is not a clique.
% \chi(G)>\omega(G)+p_0$
By \ref{bignonnbr}, there exists $v_1\in R$ such that the set of vertices in $R$ different from and nonadjacent to $v_1$
has chromatic number at least $\chi(R)/\omega(G)>m_1$, and so 
$\chi(M(v_1))>m_1$.
This proves (1).
%$\chi(G)>\omega(G)m_1+p_0$

\bigskip

Choose a stable set $S\subseteq V(G)$ with $|S|\le s$, maximal such that $\chi(N(S))>n_{|S|}$ and $\chi(M(S))>m_{|S|}$, where 
$N(S)$ denotes the set of all vertices in $V(G)\setminus S$ that are adjacent to every vertex in $S$, and $M(S)$
denotes the set of all vertices in $V(G)\setminus S$ that are nonadjacent to every vertex in $S$. From (1), $|S|\ge 1$.
Now there are two cases, $|S|<s$ and $|S|=s$.

Suppose first that $|S|<s$.
Let $A$ be the set of all vertices $v\in M(S)$ such that the set of neighbours of $v$ in $N(S)$ has chromatic number
at most $n_{|S|+1}$. Since $\chi(N(S))>n_{|S|}=qn_{|S|+1}+p$, we may assume that $\chi(A)\le p$, because otherwise
$(A,N(S))$ is a $(\psi,q)$-sprinkling. Hence 
$$\chi(B)\ge \chi(M(S))-p> m_{|S|}-p = \omega(G)m_{|S|+1},$$
where $B=M(S)\setminus A$. Since $m_{|S|+1}\ge 1$ (because $t\ge 1$), it follows that $B$ is not a clique, and so from \ref{bignonnbr},
there is a vertex $v\in B$ such that the set of vertices in $B$, different from and nonadjacent to $v$, has chromatic number
at least $\chi(B)/\omega(G)>m_{|S|+1}$. But then adding $v$ to $S$ contradicts the maximality of $S$.

Now suppose that $|S|=s$. Since $\chi(N(S))> n_s=t$, we may choose $T\subseteq N(S)$ with $|T|=t$. Let $A$ be the set of vertices in $M(S)$
that have $s$ non-neighbours in $T$ that are pairwise nonadjacent, and let $B=M(S)\setminus A$.
For each stable set $S'\subseteq T$ with $|S'|=s$,
we may assume that the set of vertices in $M(S)$ with no neighbour in $S'$ has chromatic number at most $p$, because otherwise
$G[S\cup S']$ is a $\psi$-nondominating copy of $K_{s,s}$. The number of such sets $S'$ is at most $t^s$, and so $\chi(A)\le t^sp$.
Hence 
$$\chi(B)\ge \chi(M(S))-t^sp>m_s-t^sp=t,$$ 
and so there exists $M\subseteq B$ with $|M|=t$. But then $(M,T)$ is a $(t,2)$-template.
This proves \ref{bigbip}.~\bbox

We also need the following version of Ramsey's theorem (proved for instance in~\cite{poly2}).
\begin{thm}\label{ramsey}
	For all integers $s\ge 1$ and $r\ge 2$, if a graph $G$ has no stable subset of size $s$ and no clique of size more than $r$, 
	then $|V(G)|< r^s$.
\end{thm}

Now we use \ref{bigbip} to prove \ref{strongisol}, which we restate in a strengthened form:
\begin{thm}\label{strongisol2}
        Let $s,q\ge 0$ be integers, and let $\psi:\mathbb{N}\rightarrow\mathbb{N}$ be a non-decreasing polynomial.
	Let $\phi,\phi':\mathbb{N}\rightarrow\mathbb{N}$ be the polynomials defined by
\begin{align*}
	\phi'(x)&= x^{s}\left(s\psi(x)+(s+1)^sx^{s(s+1)}\psi(x)+(s+1)x^{s+1}\right)\\
	&\ \ \ \ \ \ \ \ \ \ \ \ + q^{s} (s+1)x^{s+1}+ \left(2+q+q^2+\cdots+q^{s-1}\right)\psi(x)+2\\
	\phi(x)&=(s+1)^{2s}x^{2+2s(s+1)}\psi(x)+   (s+1)^sx^{1+s(s+1)}\phi'(x) + (x+1)(s+1)x^{s+1}.
\end{align*}
	for all $x\in \mathbb{N}$.
	Let $G$
        be a graph with $\chi(G)>\phi(\omega(G))$. Then either:
        \begin{itemize}
                \item there is a $\psi$-nondominating copy of $K_{s,s}$ in $G$, or
                \item there is a  $(\psi, q)$-sprinkling in $G$.
        \end{itemize}
\end{thm}
\Proof
Let $t=(s+1)\omega(G)^{s+1}$. Thus 
$$\chi(G)>\omega(G)^{2}t^{2s}\psi(\omega(G)+   \omega(G)t^s\phi'(\omega(G)) + (\omega(G)+1)t.$$
We claim we may assume that:
\\
\\
(1) {\em If $A\subseteq V(G)$ with $\chi(A)>\phi'(\omega(G))$ then $G[A]$ contains a $(t,2)$-template.}
\\
\\
Suppose not. Let $G'=G[A]$. Since $\chi(A)>\phi'(\omega(G))$ and $\psi$ is nondecreasing, it follows that
$$\chi(G')> \omega(G')^{s}\left(t^s\psi(\omega(G'))+t\right)+ s\omega(G')^{s}\psi(\omega(G'))$$
        and $\chi(G')\ge q^{s} t+ \left(2+q+q^2+\cdots+q^{s-1}\right)\psi(\omega(G')) +2$.
By \ref{bigbip} applied to $G'$, either
\begin{itemize}
	\item there is a $\psi$-nondominating copy of $K_{s,s}$ in $G'$ (and hence in $G$), or
	\item there is a  $(\psi, q)$-sprinkling in $G'$ (and hence in $G$), or
                \item $G'$ contains a $(t,2)$-template.
        \end{itemize}
We may assume that neither of the first two bullets hold, so the third holds. This proves (1).

\bigskip

For $2\le k\le \omega(G)+1$, define $t_k=(s+1)\omega(G)^{s+1}-s(k-2)\omega(G)^s$. Thus $t_2=t$, and 
$0\le t_k\le t$ for $2\le k\le \omega(G)+1$.
By (1) applied to $G$, there is a $(t_2,2)$-template in $G$.
Choose an integer $k$ with $2\le k\le \omega(G)+1$, maximum such that there is a $(t_k,k)$-template in $G$, and let $(A_1\LL A_k)$
be such a template.
\\
\\
(2) {\em $k\le \omega(G)$.}
\\
\\
Suppose that $k=\omega(G)+1$. Inductively for $i=1\LL k$, suppose that vertices $a_1\LL a_{i-1}$ are defined, and define $a_i$
as follows. For $1\le h<i$, the non-neighbours of $a_h$ in $A_i$ do not include a stable set of cardinality $s$, from the definition
of a $(t_k,k)$-template. Hence by \ref{ramsey} (taking $r=\omega(G)$), there are at most $\omega(G)^s$ vertices in $A_i$ nonadjacent to $a_h$, and hence
at most $\omega(G)^{s+1}$ vertices in $A_i$ that are nonadjacent to at least one of $a_1\LL a_{i-1}$. Since
$$|A_i|=t_k\ge (s+1)\omega(G)^{s+1}-s(\omega(G)-1)\omega(G)^s>\omega(G)^{s+1},$$
some vertex $a_i\in A_i$ is adjacent to all of $a_1\LL a_{i-1}$. This completes the inductive definition. But then 
$\{a_1\LL a_{\omega(G)+1}\}$ is a clique in $G$, a contradiction. This proves (2).

\bigskip

Let $Z=V(G)\setminus (A_1\cupcup A_k)$. For $1\le i\le k$, let $\mathcal{S}_i$ be the set of all stable sets contained in $A_i$
with cardinality $s$. For each $S\in \mathcal{S}_i$, let $D_S$ be the set of vertices in $Z$ with no neighbour in $S$,
and let $Y_i$ be the union of the sets $D_S$ over $S\in \mathcal{S}_i$.
\\
\\
(3) {\em $|Z\setminus (Y_1\cupcup Y_k)|< t_{k+1}$.}
\\
\\
Suppose not, and choose 
$A\subseteq Z\setminus (Y_1\cupcup Y_k)$ with $|A|=t_{k+1}$. For $1\le i\le k$, choose $B_i\subseteq A_i$ with $|B_i|=t_{k+1}$. Then 
$(A,B_1,B_2\LL B_k)$ is a $(t_{k+1}, k+1)$-template, contrary to the maximality of $k$. This proves (3).

\bigskip

For each $v\in Y_1\cupcup Y_k$, choose $i\in \{1\LL k\}$ minimum such that $v\in Y_i$, and choose $S\in \mathcal{S}_i$
such that $v\in D_S$. We call $S$ the {\em home} of $v$.
\\
\\
(4) {\em Let $1\le i\le k$, and let $S\in \mathcal{S}_i$. The set of vertices in $D_S$ with home $S$ has chromatic number at most 
$\omega(G)t^s\psi(\omega(G))+\phi'(\omega(G))$.}
\\
\\
Let $F$ be the set of vertices in $D_S$ with home $S$. By \ref{ramsey}, as in the proof of (2), for $i+1\le j\le k$ there are at most $s\omega(G)^s$
vertices in $A_j$ with a non-neighbour in $S$, and since $|A_j|=t_k=t_{k+1}+s\omega(G)^s$, there exists $B_j\subseteq A_j$ with 
$|B_j|=t_{k+1}$ complete to $S$. For $1\le h<i$, choose $B_h\subseteq A_h$ with $|B_h|=t_{k+1}$ arbitrarily. 
Let $F'$ be the set of vertices $v\in F$ such that $v$ has no neighbour in $S'$ for some $j\in \{i+1\LL k\}$ and some 
$S'\in \mathcal{S}_j$.
For $i+1\le j\le k$, and each $S'\in \mathcal{S}_j$, the chromatic number of the set of vertices in $F$ with no neighbour in $S'$
is at most $\psi(\omega(G))$, since the copy of $K_{s,s}$ induced on $S\cup S'$ is not $\psi$-nondominating; and so 
$\chi(F')\le \omega(G)t^s\psi(\omega(G))$, since there are at most $\omega(G)t^s$ choices for the pair $(j,S')$. 
Let $F''=F\setminus F'$.
If $G[F'']$
contains a $(t,2)$-template, then it contains a $(t_{k+1}, 2)$-template $(C_1,C_2)$ say; and then 
$$(C_1,C_2, B_1\LL B_{i-1}, B_{i+1}\LL B_k)$$
is a $(t_{k+1}, k+1)$-template in $G$, from the definition of a home, a contradiction. Thus $G[F'']$ contains no such template, and 
so $\chi(F'')\le \phi'(\omega(G))$ by (1). Hence $\chi(F)\le \omega(G)t^s\psi(\omega(G))+ \phi'(\omega(G))$.
This proves (4).

\bigskip

Now every vertex in $Y_1\cupcup Y_k$ has a home, and there are only at most $\omega(G)t^s$ choices of a home; so by
(4), $\chi(Y_1\cupcup Y_k)\le \omega(G)^2t^{2s}\psi(\omega(G)) + \omega(G)t^s\phi'(\omega(G))$. Hence 
\begin{align*}
\chi(G)&\le  \omega(G)^2t^{2s}\psi(\omega(G))+ \omega(G)t^s\phi'(\omega(G))+|Z\setminus (Y_1\cupcup Y_k)|+|A_1\cupcup A_k|\\
&\le 
 \omega(G)^2t^{2s}\psi(\omega(G))+ \omega(G)t^s\phi'(\omega(G)) + (\omega(G)+1)t,
\end{align*}
a contradiction. This proves \ref{strongisol2}.~\bbox

\section{Odd holes}

Now we deduce \ref{longandfourmultihole}. Let us say a
hole in $G$ is {\em special} if its length 
is either four or odd.
We need a result proved in~\cite{ind7}, the following:
\begin{thm}\label{compconj}
	Let $x\in \mathbb{N}$, and let $G$ be a
        graph such that $\chi(N(v))\le x$ for every vertex $v\in V(G)$. 
	If $C$ is a shortest odd hole in $G$, the set of vertices of $G$ that belong to or have a neighbour in $V(C)$
	has chromatic number at most $21x$.
\end{thm}
We deduce:
\begin{thm}\label{smallnbrs}
        Let $\psi:\mathbb{N}\rightarrow\mathbb{N}$ be some non-decreasing polynomial, let $n\in \mathbb{N}$, and let $G$ be a
	graph such that $\chi(N(v))\le n$ for every vertex $v\in V(G)$. If $\chi(G)> \max(\omega(G),21n+\psi(\omega(G)))$ 
	then $G$ contains a 
	$\psi$-nondominating special hole.
\end{thm}
\Proof
Since $\chi(G)>\omega(G)$, $G$ is not perfect, and so contains either  a four-hole or an odd hole (by the strong perfect graph theorem~\cite{SPGT}, since odd antiholes of length at least seven contain four-holes). Let $C$
be either a four-hole, or a shortest odd hole of $G$. Let $A$ be the set of vertices in $V(G)\setminus V(C)$ that have no
neighbour in $V(C)$, and $B=V(G)\setminus A$. If $C$ has length four then $\chi(B)\le 4n$, and if $C$ is 
a shortest odd hole of $G$, then $\chi(B)\le 21n$ by \ref{compconj}. Consequently 
$\chi (A)>\psi(\omega(G))\ge \psi(\omega(A))$,
and so $C$ is a $\psi$-nondominating special hole. This proves \ref{smallnbrs}.~\bbox

We also need:
\begin{thm}\label{bignbrs}
	Let $G$ be a
	graph containing no four-hole, let $n\in \mathbb{N}$,  and let $X\subseteq V(G)$ be the set of all $v\in V(G)$ with 
	$\chi(N(v))> n$. If $\chi(X)> \omega(G)$, then
	there exist disjoint sets $A,B\subseteq V(G)$, anticomplete, with $\chi(A), \chi(B)>n/2-\omega(G)$.
\end{thm}
\Proof
Let us say an edge $xy$ of $G$ is {\em rich} if $\chi(N(x) \setminus N(y))>n/2 -\omega(G)$ and 
$\chi(N(y) \setminus N(x))>n/2 -\omega(G)$.
Since there is no four-hole, it is enough to prove that there is a rich edge.

Since $\chi(X)> \omega(G)$, the graph $G[X]$ is not perfect, and so contains a four-vertex induced path with vertices
$v_1\DD v_2 \DD v_3 \DD v_4$ in order.
Let 
\begin{align*}
    A_1=&N(v_1) \setminus (N(v_3) \cup N(v_4))\\
    A_2=&N(v_2) \setminus (N(v_4) \cup (N(v_1)\cap N(v_3)))\\
    A_3=&N(v_3) \setminus (N(v_1) \cup (N(v_2)\cap N(v_4)))\\
    A_4=& N(v_4) \setminus (N(v_2) \cup N(v_1)).
\end{align*}

Since there is no four-hole, $N(v_1)\cap N(v_3)$ is a clique, and so is $N(v_1)\cap N(v_4)$, and therefore $\chi(A_1) > n- 2\omega(G)$.
Since $N(v_2) \cap N(v_4)$ and $N(v_1)\cap N(v_3))$ are cliques, it also follows that $\chi(A_2) > n- 2\omega(G)$, and 
similarly $\chi(A_i) > n- 2\omega(G)$ for $1\le i\le 4$.

Now $v_2$ is anticomplete to $A_1\setminus A_2$, and $v_1$ is anticomplete to $A_2\setminus A_1$, so
if $\chi(A_1 \cap A_2) \le  n/2 - \omega(G)$, then $\chi(A_1\setminus A_2) > n/2-\omega(G)$ and
$\chi(A_2 \setminus  A_1) > n/2 - \omega(G)$, and so the edge $v_1v_2$ is rich.

Thus we may assume that $\chi(A_1 \cap A_2) > n/2 - \omega(G)$, and similarly $\chi (A_3 \cap A_4) > n/2-\omega(G)$.
But  $A_1 \cap A_2 \subseteq N(v_2) \setminus N(v_3)$, and
$A_3 \cap A_4 \subseteq  N(v_3) \setminus N(v_2)$,
and so the edge $v_2v_3$ is rich. This proves \ref{bignbrs}.~\bbox

We put \ref{smallnbrs} and \ref{bignbrs} together to prove the following:
\begin{thm}\label{girth5}
        Let $\psi:\mathbb{N}\rightarrow\mathbb{N}$ be some non-decreasing polynomial. If $G$ is a $C_4$-free graph 
	with
	$$\chi(G)>85\omega(G)+43\psi(\omega(G))$$ 
	then $G$ contains 
	a $\psi$-nondominating odd hole.
\end{thm}
\Proof
Let $G$ be a $C_4$-free graph with 
%$c\ge 1$
$\chi(G)>85\omega(G)+43\psi(\omega(G))$.
Define $n=4\omega(G)+2\psi(\omega(G))$.

Let $A$ be the set of all vertices $v$ of $G$ such that $\chi(N(v))\le n$, and $B=V(G)\setminus A$. By \ref{smallnbrs}
applied to $G[A]$, we may assume that 
$$\chi(A)\le \max(\omega(A),21n+\psi(\omega(A)))=21n+\psi(\omega(A))\le 84\omega(G)+43\psi(\omega(G))$$
and so $\chi(B)\ge \chi(G)-\chi(A)>\omega(G)$.
By \ref{bignbrs} 
there exist disjoint sets $X,Y\subseteq V(G)$, anticomplete, with $\chi(X), \chi(Y)>n/2-\omega(G)\ge \omega(G)+\psi(\omega(G))$.
Since $\chi(X)>\omega(G)\ge \omega(X)$, $G[X]$ is not perfect and so contains a special hole $C$, and hence an odd hole since
$G$ has no four-holes. Since $V(C)$ is anticomplete to $Y$,
and $\chi(Y)>\psi(\omega(G))\ge \psi(\omega(Y))$, $C$ is $\psi$-nondominating. This proves \ref{girth5}.~\bbox

This in turn is used to prove:

\begin{thm}\label{nondomhole} Let $\psi:\mathbb{N}\rightarrow\mathbb{N}$ be some non-decreasing polynomial. Then there is
	a non-decreasing polynomial $\phi:\mathbb{N}\rightarrow\mathbb{N}$
	such that if $\chi(G)>\phi(\omega(G))$ then $G$ contains a $\psi$-nondominating special hole.
\end{thm}
\Proof
Let $\psi'(x)=85x+43\psi(x)$ for $x\in \mathbb{N}$, and 
let $\phi$ satisfy \ref{strongisol} with $\psi$ replaced by $\psi'$, taking $s=2$ and $q=4$. 
We claim that $\phi$ satisfies \ref{nondomhole}.
Thus, let 
$G$ be a graph with $\chi(G)>\phi(\omega(G))$.  By \ref{strongisol}, either
there is a $\psi'$-nondominating four-hole in $G$, or
there is a  $(\psi', 4)$-sprinkling in $G$. In the first case, this four-hole is also $\psi$-nondominating, since $\psi(x)\le \psi'(x)$
for $x\in \mathbb{N}$, so we assume the second case holds. Let $(P,Q)$ be a $(\psi', 4)$-sprinkling in $G$, and let $r$ be the
maximum chromatic number over $v\in P$ of the set of neighbours of $v$ in $Q$. Thus $\chi(Q)> 4r+\psi'(\omega(Q))$,
from the definition of a $(\psi', 4)$-sprinkling. If $G[P]$ has a four-hole $H$,
the set of vertices in $Q$ with a neighbour in $V(H)$ has chromatic number at most $4r$, and so there is a subset of $Q$ 
with chromatic number more than $\psi'(\omega(Q))\ge \psi(\omega(Q))$ anticomplete to $H$, and so $H$ is $\psi$-nondominating.
Thus we may assume that $G[P]$ has no four-hole. By \ref{girth5}, $G[P]$, and hence $G$, contains a $\psi$-nondominating odd hole.
This proves \ref{nondomhole}.~\bbox

We deduce \ref{oddandfourmultihole}, which we restate:
\begin{thm}\label{oddandfourmultihole2}
        For each integer $k\ge 0$, let $\mathcal{C}$ be the class of all graphs $G$ with
        no $k$-multihole in which every component is special. Then $\mathcal{C}$ is \poly.
\end{thm}
\Proof
Let us say a $k$-multihole is {\em special} if each of its components is a special hole. We proceed by induction on $k$. 
The result is true when $k=1$, because graphs containing no special hole are perfect;
so we assume that $k\ge 2$, and there is a polynomial  binding function
$\psi:\mathbb{N}\rightarrow\mathbb{N}$ for the class of all graphs with no special $(k-1)$-multihole
$\mathcal{C}_{k-1}$ (and we may assume $\psi$ is non-decreasing).
Let $\phi$ satisfy \ref{nondomhole}; 
we claim that $\phi$ is a binding function for the class of all graphs with no special $k$-multihole.
Thus, let $G$ be a graph with $\chi(G)>\phi(\omega(G))$; we must show that $G$ contains a special $k$-multihole. By 
the choice of $\phi$, $G$ contains a $\psi$-nondominating special hole $H$ say. Choose $A\subseteq V(G)\setminus V(H)$,
anticomplete to $V(H)$, such that $\chi(A)> \psi(\omega(A))$. From the inductive hypothesis, $G[A]$ contains a special $(k-1)$-multihole,
and so $G$ contains a special $k$-multihole. This proves \ref{oddandfourmultihole2}.~\bbox

\section{Long holes}

In this section we will prove \ref{longandbip}. The proof is similar to that of \ref{oddandfourmultihole}.
Fix an integer $\ell\ge 4$, and we say a hole is {\em long} if its length is at 
least $\ell$.
Let $\tau(G)$ denote the largest integer $t$ such that $G$ contains $K_{t,t}$ as a subgraph.
We need a result proved in~\cite{bonamy} (see also \cite{poly1}), the following:

\begin{thm}\label{longhole}
	There exists an integer $c>0$ such that $\chi(G)\le \tau(G)^c+1$ for
every graph $G$ with no long hole.
\end{thm}
We deduce:
\begin{thm}\label{biblonghole}
	Let $s\in \mathbb{N}$; then the class of $K_{s,s}$-free graphs with no long hole is \poly.
\end{thm}
\Proof
Let $c\ge 1$ be as in \ref{longhole}, and let $\phi$ be the polynomial $\phi(x)=x^{cs}$ for $x\in \mathbb{N}$.
Let $G$ be a  $K_{s,s}$-free graph with no long hole. We will show that $\phi$ is a binding function for $G$.
Suppose that $\tau(G)\ge \omega(G)^s$, and let $A,B$ be disjoint subsets of $V(G)$,
both of cardinality at least $\omega(G)^s$ and complete to each other. By \ref{ramsey}, there exist stable sets $A'\subseteq A$
and $B'\subseteq B$ both of cardinality $s$; but then $G[A'\cup B']$ is a copy of $K_{s,s}$, a contradiction. So  $\tau(G)< \omega(G)^s$.
By \ref{longhole}, 
$$\chi(G)\le (\omega(G)^s-1)^c+1\le \omega(G)^{cs}=\phi(\omega(G)),$$ 
and so $\phi$ is a binding function for $G$, and hence for 
the class of $K_{s,s}$-free graphs with no long hole. This proves \ref{biblonghole}.~\bbox

Next we need an analogue of \ref{smallnbrs}, the following:
\begin{thm}\label{smalllong}
        Let $n\in \mathbb{N}$, and let $G$ be a
        graph such that $\chi(N(v))\le n$ for every vertex $v\in V(G)$.
        If $C$ is a shortest long hole in $G$, the set of vertices of $G$ that belong to or have a neighbour in $V(C)$
	has chromatic number at most $(\ell+1)n$.
\end{thm}
\Proof
Let $C$ have vertices $c_1\DD c_2 \CC c_k\DD c_1$ in order. Let $P$ be the path $c_1\DD c_2 \CC c_{\ell-3}$, and let 
$Q$ be the path $C\setminus V(P)$. 
\\
\\
(1) {\em If $v\in V(G)\setminus V(C)$ has no neighbour in $V(P)$, then 
all neighbours of $v$ in $V(Q)$ belong to a three-vertex subpath of $Q$.}
\\
\\
Suppose not, and choose $i,j$ minimum and maximum respectively such that $c_i,c_j\in V(Q)$ are neighbours of $v$. Thus $j-i\ge 3$,
and so 
$$c_1\DD c_2\CC c_i\DD v\DD c_j\DD c_{j+1}\CC c_k\DD c_1$$
is a long hole (because $j\ge \ell-2$) that is shorter than $C$, a contradiction. This proves (1).

\bigskip

For $1\le i\le k$, let $A_i$ be the set of vertices in $V(G)\setminus V(C)$ that are adjacent to $c_i$ and to none of
$c_1\LL c_{i-1}$.
\\
\\
(2) {\em $A_i$ is anticomplete to $A_j$ for $\ell-2\le i<j\le k$ with $j-i\ge 4$.}
\\
\\
Suppose that $u\in A_i$ and $v\in A_j$ are adjacent. Choose $j'\ge j$ maximum such that $c_{j'}$ is adjacent to $v$; thus 
$j' \ge j \ge i+4$, and so by (1), $u$ is
non-adjacent to $c_{j'}\LL c_k$. Hence
$$c_1\DD c_2\CC c_i\DD u\DD v\DD c_{j'}\DD c_{j'+1}\CC c_k\DD c_1$$
is a long hole shorter than $C$, a contradiction. This proves (2).

\bigskip
For $t = 1,2,3,4$ let $I_t$ be the set of all integers $i\in \{\ell-2\LL k\}$ such that $i-t$ is divisible by four. Thus
$I_1,I_2,I_3,I_4$ form a partition of $\{\ell-2\LL k\}$. Moreover, for all $t\in \{1\LL 4\}$, and all distinct $i,j\in I_t$,
there is no edge between $A_i\cup \{c_{i+1}\}$ and $A_j\cup \{c_{j+1}\}$, by (2); and so 
$\bigcup_{i\in I_t}A_i\cup \{c_{i+1}\}$ has chromatic number at most $n$. Hence the set of all vertices in $V(G)$ that belong to or have
a neighbour in $V(C)$ has chromatic number at most $(\ell+1)n$, since those that belong to or have a neighbour in $P$
have chromatic number at most $(\ell-3)n$, and the others have chromatic number at most $4n$. This proves \ref{smalllong}.~\bbox

Now we need an analogue of \ref{bignbrs}, the following:
\begin{thm}\label{biglong}
	Let $s\in \mathbb{N}$, let $G$ be a $K_{s,s}$-free graph, with no long hole of length at most $2s\ell$.
        Let $n\in \mathbb{N}$,  and let $B\subseteq V(G)$ be the set of vertices $v$
	of $G$ such that $\chi(N(v))> n$. If $G[B]$ contains a long hole, then
		there exist disjoint sets $X,Y\subseteq B$, anticomplete, with $\chi(X), \chi(Y)>n-(2s\ell)^s\omega(G)^s$.
\end{thm}
\Proof We may assume that $G[B]$ has a hole of length more than $2s\ell$, and so contains an induced path $P$ with $2s\ell-1$
vertices. Let the vertices of $P$ be $p_1\DD p_2\CC p_{r}$ in order, where $r=2s\ell-1$. For each stable subset 
$S\subseteq V(P)$ with $|S|=s$, let $D_S$
be the set of vertices in $V(G)\setminus V(P)$ that are adjacent to every vertex in $S$. Since $G$ is $K_{s,s}$-free, it follows 
from \ref{ramsey} that $|D_S|\le \omega(G)^s$. Let $D$ be the set of vertices in 
$V(G)\setminus V(P)$ that have $s$ pairwise nonadjacent neighbours in $V(P)$. Since there are at most $(2s\ell)^s$ choices of $S$, it
follows that $\chi(D)\le (2s\ell)^s\omega(G)^s$.
Let $F=V(G)\setminus (V(P)\cup D)$. 
\\
\\
(1) {\em For each $v\in F$, if $i,j$ are minimum and maximum such that $v$ is adjacent to $p_i, p_j$, then $j-i\le (s-2)(\ell-2)+1$.}
\\
\\
Let $v\in F$. Choose $t\ge 0$ maximum such that there exist $1\le i_1<\cdots<i_t\le r$ satisfying:
\begin{itemize}
	\item $i_1$ is the least $i$ such that $v$ is adjacent to $p_i$;
	\item $v$ is adjacent to $p_{i_k}$ for $1\le k\le t$;
	\item $i_{k+1}\ge i_k+2$ for $1\le k\le t-1$; 
	\item $v$ is nonadjacent to $p_j$ for $1\le k\le t-1$ and for each $j\in \{i_k+2\LL i_{k+1}-1\}$.
\end{itemize}
Since $\{p_{i_1},p_{i_2}\LL p_{i_t}\}$ is a stable set, and $v\in F$,
it follows that $t<s$. Moreover, for $1\le k<t$, $v$ is nonadjacent to each $p_j$ for each $j\in \{i_k+2\LL i_{k+1}-1\}$; so one
of 
$$v\DD p_{i_k}\DD p_{i_k+1}\CC p_{i_{k+1}}$$
$$v\DD p_{i_k+1}\DD p_{i_k+2}\CC p_{i_{k+1}}$$
is an induced cycle. This cycle has length at most $2s\ell$, since $P$ has only $r=2s\ell-1$ vertices; and so the cycle 
has length
less than $\ell$, since $G$ has no long hole of length at most $2s\ell$. Consequently $i_{k+1}-i_k\le \ell-2$, and so $i_t-i_1\le (s-2)(\ell-2)$.
From the maximality of $t$, $v$ is nonadjacent to $p_j$ for all $j\ge i_t+2$. This proves (1).

\bigskip

Let $X$ be the set of neighbours of $p_1$ in $V(G)\setminus D$, and let $Y$ be the set of neighbours of $p_{r}$ in $V(G)\setminus D$.
\\
\\
(2) {\em $X$ is disjoint from and anticomplete to $Y$.}
\\
\\
Since $r-1>(s-2)(\ell-2)+1$, (1) implies that $X\cap Y=\emptyset$. Suppose that $u\in X$ and $v\in Y$ are adjacent. 
Choose $i\in \{1\LL r\}$ maximum such that $u$ is adjacent to 
$p_i$, and choose $j\in \{1\LL r\}$ minimum such that $v$ is adjacent to
$p_j$. By (1), $i-1\le (s-2)(\ell-2)+1$, and $r-j\le (s-2)(\ell-2)+1$. Hence 
$i-1+r-j\le 2((s-2)(\ell-2)+1)$, and so 
$$j-i\ge (r-1)-2((s-2)(\ell-2)+1)=4\ell+4s-12.$$
But then $u\DD p_i\DD p_{i+1}\CC p_j\DD v\DD u$ is a hole of length at least $4\ell+4s-9\ge \ell$ and at most $2s\ell$,
a contradiction. This proves (2).

\bigskip

But $\chi(N(p_1))\ge n$, and so $\chi(X)\ge n-\chi(D)\ge n-(2s\ell)^s\omega(G)^s$, and the same for $Y$. This proves \ref{biglong}.~\bbox

Next, combining \ref{smalllong} and \ref{biglong}, we have an analogue of \ref{girth5}:

\begin{thm}\label{longgirth5}
	Let $s\in \mathbb{N}$, and let $\psi:\mathbb{N}\rightarrow\mathbb{N}$ be some non-decreasing polynomial. There is a 
	non-decreasing polynomial $\phi:\mathbb{N}\rightarrow\mathbb{N}$ with the following property.
	If $G$ is a $K_{s,s}$-free
	graph with no long hole of length at most $2s\ell$, and no $\psi$-nondominating long hole, then $\chi(G)\le \phi(\omega(G))$.
\end{thm}
\Proof
By \ref{biblonghole}, there is a non-decreasing polynomial $\theta:\mathbb{N}\rightarrow\mathbb{N}$ that is a binding function
for the class of $K_{s,s}$-free graphs with no long hole. Define $\phi$ by 
$$\phi(x)=2\theta(x)+\psi(x)+(\ell+1)\left((2s\ell)^sx^s+\theta(x)+\psi(x)\right).$$
We claim that $\phi$ satisfies \ref{longgirth5}. Thus, let $G$ be  a $K_{s,s}$-free
graph with no long hole of length at most $2s\ell$, and no $\psi$-nondominating long hole. Let 
$$n=(2s\ell)^s\omega(G)^s+\theta(\omega(G))+\psi(\omega(G)).$$ 
Let $A$ be the set of vertices $v\in V(G)$ such that $\chi(N(v))\le n$, and $B=V(G)\setminus A$. 
\\
\\
(1) {\em $\chi(A)\le \theta(\omega(G))+\psi(\omega(G))+(\ell+1)n$.}
\\
\\
Suppose not. Then by \ref{biblonghole}, $G[A]$ has a long hole; let $C$ be a shortest long hole of $G[A]$. By \ref{smalllong}
applied to $G[A]$,
the set of vertices of $A$ that belong to or have a neighbour in $V(C)$
	has chromatic number at most $(\ell+1)n$, and so there is a subset of $A\setminus V(C)$ anticomplete to $V(C)$
	with chromatic number more than $\chi(A)- (\ell+1)n\ge \psi(\omega(G))$. Hence $C$ is $\psi$-nondominating, 
	a contradiction. This proves (1).
\\
\\
(2) {\em $\chi(B)\le  \theta(\omega(G)) $.}
\\
\\
Suppose not. Then $G[B]$ has a long hole by \ref{biblonghole}. By \ref{biglong},
there exist disjoint sets $X,Y\subseteq B$, anticomplete, with $\chi(X), \chi(Y)>n-(2s\ell)^s\omega(G)^s$.
Since $\chi(X)\ge \theta(\omega(G))$, $G[X]$ has a long hole, and it is $\psi$-nondominating since $\chi(Y)\ge \psi(\omega(G))$,
a contradiction.
%x-(2s\ell)^s\omega(G)^s\ge q_1(\omega(G))$
%x-(2s\ell)^s\omega(G)^s\ge p(\omega(G))$
This proves (2).

\bigskip

From (1) and (2), it follows that 
$$\chi(G)\le 2\theta(\omega(G))+\psi(\omega(G))+(\ell+1)n.$$ 
This proves \ref{longgirth5}.~\bbox

This implies:

\begin{thm}\label{nondomlonghole} 
	Let $s\in \mathbb{N}$, and let $\psi:\mathbb{N}\rightarrow\mathbb{N}$ be some non-decreasing polynomial. Then there is
        a non-decreasing polynomial $\phi:\mathbb{N}\rightarrow\mathbb{N}$
        such that if $\chi(G)>\phi(\omega(G))$ then $G$ contains either  a $\psi$-nondominating
	copy of $K_{s,s}$, or a $\psi$-nondominating long hole.
\end{thm}
\Proof
By \ref{longgirth5},
there is a
        non-decreasing polynomial $\psi':\mathbb{N}\rightarrow\mathbb{N}$ with the following property.
        If $G$ is a $K_{s,s}$-free
        graph with no long hole of length at most $2s\ell$, and $\chi(G)>\psi'(\omega(G))$,
        then $G$ contains
        a $\psi$-nondominating long hole.

Let $\phi$ satisfy \ref{strongisol} with $\psi$ replaced by $\psi'$, taking $q=2s\ell$.
We claim that $\phi$ satisfies \ref{nondomlonghole}.
Thus, let
$G$ be a graph with $\chi(G)>\phi(\omega(G))$.  By \ref{strongisol}, either
there is a $\psi'$-nondominating copy of $K_{s,s}$ in $G$, or
there is a  $(\psi', 2s\ell)$-sprinkling in $G$. In the first case, this copy of $K_{s,s}$ is also $\psi$-nondominating, since $\psi(x)\le \psi'(x)$
for $x\in \mathbb{N}$, so we assume the second case holds. Let $(P,Q)$ be a $(\psi', 2s\ell)$-sprinkling in $G$, and let $r$ be the
maximum chromatic number over $v\in P$ of the set of neighbours of $v$ in $Q$. Thus $\chi(Q)> 2s\ell r+\psi'(\omega(Q))$,
from the definition of a $(\psi', 2s\ell)$-sprinkling. If $G[P]$ contains $H$ where $H$ is either a copy of $K_{s,s}$ or a long hole of length at most $2s\ell$,
the set of vertices in $Q$ with a neighbour in $V(H)$ has chromatic number at most $|H|r\le 2s\ell r$, and so there is a subset of $Q$
with chromatic number more than $\psi'(\omega(Q))\ge \psi(\omega(Q))$ anticomplete to $H$; and therefore $H$ is $\psi$-nondominating.
Thus we may assume that $G[P]$ is $K_{s,s}$-free and has no long hole of length at most $2s\ell$. By \ref{longgirth5}, $G[P]$, and hence $G$, contains a $\psi$-nondominating long hole.
This proves \ref{nondomlonghole}.~\bbox

Finally, we prove \ref{longandbip}, which we restate:
\begin{thm}\label{longandbip2}
        For all integers $k, s\ge 0$ and $\ell\ge 4$, let $\mathcal{C}$ be the class of all graphs $G$ such that no induced subgraph of $G$
        has exactly $k$ components, each of which is either a copy of $K_{s,s}$ or a cycle of length at least $\ell$.
        Then $\mathcal{C}$ is \poly.
\end{thm}
\Proof
(The proof is just like that of \ref{oddandfourmultihole2}.)
Let us say an induced subgraph $H$ of a graph $G$ is a {\em $k$-object} if it has exactly $k$ components, and each is either
a copy of $K_{s,s}$ or a cycle of length at least $\ell$. Thus $\mathcal{C}_k$ is the class of graphs with no $k$-object.
We prove by induction on $k$ that $\mathcal{C}_k$ is \poly.
The result is true when $k=1$, by \ref{biblonghole},
so we assume that $k\ge 2$, and there is a polynomial binding function
$\psi:\mathbb{N}\rightarrow\mathbb{N}$ for 
$\mathcal{C}_{k-1}$ (and we may assume $\psi$ is non-decreasing).
Let $\phi$ satisfy \ref{nondomlonghole};
we claim that $\phi$ is a binding function for $\mathcal{C}_k$.
Thus, let $G$ be a graph with $\chi(G)>\phi(\omega(G))$; we must show that $G$ contains a $k$-object. By
the choice of $c$, $G$ contains a $\psi$-nondominating induced subgraph $H$, where $H$ is either a copy of $K_{s,s}$ or a long hole.
Choose $A\subseteq V(G)\setminus V(H)$,
anticomplete to $V(H)$, such that $\chi(A)> \psi(\omega(A))$. From the inductive hypothesis, $G[A]$ contains a $(k-1)$-object,
and so $G$ contains a $k$-object. This proves \ref{longandbip2}.~\bbox

%RR

\end{document}